\documentclass[10pt]{article}
\usepackage{amsmath}
\usepackage{amssymb,latexsym}
\usepackage{amsfonts}
\usepackage{times}
\usepackage{natbib}
\usepackage{indentfirst}
\usepackage[dvips]{graphicx}
\newtheorem{proposition}{Proposition}
\usepackage{fullpage}
\usepackage{setspace}
\title{Estimating abundance-based generalized species accumulation curves}
\date{}
\author{Chang Xuan Mao}

\begin{document}

\maketitle

\begin{center}
Department of Statistics, University of California\\
Riverside, CA, 92521 USA\\
cmao@statserv.ucr.edu
\end{center}

\begin{abstract}
The number of species can be estimated by sampling individuals
from a species assemblage. The problem
of estimating generalized species accumulation curve 
is addressed in a nonparametric Poisson mixture model.
A likelihood-based estimator is proposed
and illustrated by real examples.
\end{abstract}


\section{Introduction}

An important but difficult problem in ecological studies
is estimating species richness, i.e., the number of species
in an assemblage based on an incomplete survey (Colwell and
Coddington 1994). The same problem also 
arises from various other scientific fields (Bunge and Fitzpatrick 1993).
In the survey, individuals are selected from the species assemblage 
and their species identities are recognized. 
The species accumulation curve (SAC) is the plot of the expected
number of species against the measure of sampling effort,
which serves a variety of purposes in ecological studies 
such as comparison among species assemblages and prediction
of expected number of new species 
(e.g., Hurlbert 1971; Colwell and Coddington 1994;
Shen et al. 2003; Mao 2005). The estimand of a nonparametric species
richness estimator is also often plotted against the measure of sampling effort,
called a generalized SAC and 
and used like the usual SAC (Colwell and Coddington 1994). 
Although estimating the usual SAC has been extensively studied
(e.g., Mao 2005), little investigation has been made to estimate generalized SACs.
A computationally intensive randomization procedure is usually used
by ecologists and conservation biologists.

Consider a species assemblage consisting of $s$ distinct
species labeled by $i=1$, 2, \dots, $s$.
The sampling of individuals from species $i$ is often modeled as 
a Poisson process with rate $\lambda_i$ over time $t\in[0,\infty)$ 
(e.g., Efron and Thisted 1976; Norris and Pollock 1998; Mao 2004, 2005).
Let $Y_i(t)$ be the number of individuals from species $i$ during  $[0,t]$. 
Conditioning on $h(t)=\sum_{i=1}^sY_i(t)$, 
the $Y_i(t)$ arise as a multinomial sample of size $h(t)$ 
with index $s$ and probabilities $p_i=\lambda_i/\sum_{j=1}^s\lambda_j$ (e.g., Chao 1984). 
When the rates $\lambda_i$ are assumed to arise as a random sample from a mixing 
distribution $\Theta=\sum_{u=1}^{\nu}\pi_u\delta(\gamma_u)$,
where $\delta(\lambda)$ is a distribution degenerate at $\lambda$,
the $Y_i(t)$ become a random sample from a Poisson mixture (e.g., Mao 2004). 

Let $n_j(t)=\sum_{i=1}^sI(Y_i(t)=j)$, where $I(\cdot)$ is the indicator function. 
Let $n(t)=(n_1(t),n_2(t),\dots)$ 
and $\phi(t)=E\{n(t)\}=(\phi_1(t),\phi_2(t),\dots)$, where
\begin{equation}\label{eqn:ecount}
\phi_j(t)=E\{n_j(t)\}=s\sum_{u=1}^{\nu}\pi_u\exp(-\gamma_ut)(\gamma_ut)^j(j!)^{-1}.
\end{equation}
Let $n_+(t)$ be the number of observed species with expectation $\phi_+(t)$, where
\[
n_+(t)=\sum_{j=1}^{\infty}n_j(t), 
\phi_+(t)=\sum_{j=1}^{\infty}\phi_j(t).
\]
A nonparametric estimator for the number of species $s$ is
a function $G(n(t))$ which
estimates $G(\phi(t))$, a parameter that approximates  $s$.
Note that $n_+(t)$ is such an estimator.
Another example is the estimator in Chao (1984),
\[
G_c(n(t))=\sum_{j=1}^{\infty}n_j(t)+\frac{n_1^2(t)}{2n_2(t)}.
\]

When the sampling is stopped at $t=t_0$, one has
a vector of observed counts $n(t_0)$. 
We will consider
the problem of estimating $G(\phi(t))$ based on $n(t_0)$.
The special case of estimating $\phi_+(t)$
was considered by Good and Toulmin (1956), Efron and Thisted (1976), 
Shen et al. (2003) and Mao (2005).

The problem can be reduced to estimating $\phi(t)$.
Good and Toulmin (1956) provided an estimator for $\phi(t)$.
The Good-Toulmin estimator usually behaves badly at $t>2t_0$
and often produces inadmissible values (e.g., negative values)
for $t\in(t_0,2t_0]$. We will develop a likelihood-based estimator,
which competes with the Good-Toulmin estimator at $t\in[0,2t_0]$
as its smoothed version. The likelihood-based estimator
is particularly useful when the Good-Toulmin estimator fails. 
Our approach is different from that in Norris and Pollock (1998)
because we do not require an estimator for $s$,
a parameter that is difficult to estimate.
We will also show that the commonly used randomization procedure is unnecessary
because it is a simulation-based approximation to an enumeration procedure
which yields an  estimator close to the Good-Toulmin estimator.

The estimation methods are detailed in Section 2. Numeric
studies are reported in Section 3. The proofs are provided 
in the Appendix. The R codes are available from the author on request.

\section{Methods}
For notational convenience, we will assume that time is scaled such that $t_0=1$.
Therefore, the full likelihood $p_0(s,\Theta)$ is given by 
\begin{equation*}
p_0(s,\Theta)=\frac{s!}{\{s-n_+(1)\}!\prod_{j=1}^{\infty}n_{j}(1)!}
g_{\Theta}^{s-n_+(1)}(0)\prod_{j=1}^{\infty} g_{\Theta}^{n_{j}(1)}(j),
\end{equation*}
where $g_{\Theta}$ is a mixture of Poisson densities,
\[
g_{\Theta}(j)=\sum_{u=1}^{\nu}\pi_u \exp(-\gamma_u)\gamma_u^j(j!)^{-1} , j=0, 1,\dots.
\]

The Good-Toulmin estimator $\tilde{\phi}_j(t)$ can be written as
\begin{equation}\label{eqn:GT}
\tilde{\phi}_j(t)=\sum_{k=0}^{\infty}\binom{k+j}{j}t^j(1-t)^k n_{k+j}(1).
\end{equation}
This estimator can arise from the following identity
\begin{equation}\label{eqn:GTsource}
\phi_j(t)=\sum_{k=0}^{\infty}\binom{k+j}{j}t^j(1-t)^k \phi_{k+j}(1),
\end{equation}
when one estimate $\phi_{x}(1)=sg_{\Theta}(x)$ by $n_x(1)$.

Let  $d=\max\{j:n_j(1)>0\}$. We can write $\tilde{\phi}_j(t)$ as
\begin{equation}\label{eqn:rep1}
\tilde{\phi}_j(t)=\sum_{b=j}^{d}\binom{b}{j}t^j(1-t)^{b-j} n_{b}(1).
\end{equation}
The last term of the series in (\ref{eqn:rep1}) dominates soon after $t>2$,
and $\tilde{\phi}_j(t)$ diverges to infinity or minus infinity
as $t$ increases, depending on whether $d-j$ is even or odd.
This might invite one to replace both $s$ and $\Theta$ with their
estimators in $\phi_j(t)$. For example, Norris and Pollock (1998)
provided nonparametric likelihood estimators for $s$ and $\Theta$
by a procedure that is computationally very expensive.

Because $s$ is difficult to estimate
(e.g., Bunge and Fitzpatrick 1993),
we will show that 
estimating $\phi(t)$ does not necessarily require
an estimator for $s$.
Note that  $p_0(s,\Theta)=p_1(s,\Theta)p_2(\Theta,n_+(1))$,
where $p_1(s,\Theta)$ is the binomial density of $n_+(1)$ and
$p_2(\Theta,n_+(1))$ is the multinomial density of $n(1)$ given $n_+(1)$, 
\begin{align*}
p_1(s,\Theta)&=\frac{s!}{\{s-n_+(1)\}!n_+(1)!}
g_{\Theta}^{s-n_+(1)}(0) \{ 1-g_{\Theta}(0) \}^{n_+(1)},\\
p_2(\Theta,n_+(1))&=\frac{n_+(1)!}{\prod_{j=1}^{\infty}n_j(1)!}
\prod_{j=1}^{\infty} \left\{\frac{g_{\Theta}(j)}{1-g_{\Theta}(0)}\right\}^{n_j(1)}.
\end{align*}
We will reformulate $p_2(\Theta,n_+(1))$ by introducing
$Q=\sum_{u=1}^{\nu}\omega_u\delta(\gamma_u)$, where
\begin{equation*}
\omega_u=\frac{\pi_u\{1-\exp(-\gamma_u)\}}
{\sum_{w=1}^{\nu}\pi_w\{1-\exp(-\gamma_w)\}}.
\end{equation*}
Let $f_Q$ be a mixture of zero-truncated Poisson densities, where
\[
f_Q(j)=\sum_{u=1}^{\nu} \omega_u \frac{\gamma_u^j}{\{\exp(\gamma_u)-1\}j!}, j\geq 1.
\]
Because it can be shown that $f_Q(j)=g_{\Theta}(j)/\{1-g_{\Theta}(0)\}$
(e.g., Mao 2004), we can rewrite $p_2(\Theta,n_+(1) )$ as $L(Q,n_+(1))$, where
\begin{equation*}
L(Q,n_+(1))=\frac{n_+(1)!}{\prod_{j=1}^{\infty}n_j(1)!}\prod_{j=1}^{\infty} f_Q^{n_j(1)}(j).
\end{equation*}

\begin{proposition}
For $j=1$, 2, \dots, $h$, and $h=1$, 2, \dots, 
\begin{equation}\label{eqn:rep2}
\phi_j(t)=\phi_+(1)\theta_j(t,Q),
\end{equation}
where $\theta_j(t,Q)$ a functional of the mixing distribution $Q$,
\begin{equation*}
\theta_j(t,Q)=\sum_{u=1}^{\nu} 
\omega_u\frac{\exp(-\gamma_u t)(\gamma_ut)^j}{\{1-\exp(-\gamma_u)\}j!}
\end{equation*}
\end{proposition}

The nonparametric maximum likelihood estimator (NPMLE) denoted by
$\widehat{Q}=\sum_{u=1}^{\hat{\nu}}\hat{\omega}_u\delta(\hat{\gamma}_u)$
maximizes $L(Q,n_+(1))$ (Lindsay 1983; Mao 2004). Because $n_+(1)$ estimates
$\phi_+(1)$, from (\ref{eqn:rep2}),
a likelihood-based estimator $\hat{\phi}_j(t)$ for $\phi_j(t)$ is given by 
\begin{align}
\hat{\phi}_j(t)=n_+(1)\theta_j(t,\widehat{Q}).\label{eqn:est2}
\end{align}
Note that $\hat{\phi}_j(t)$ is a smoothed version of 
$\tilde{\phi}_j(t)$ in (\ref{eqn:GT}) because
\begin{equation}\label{eqn:GT2}
\hat{\phi}_j(t)=\sum_{k=0}^{\infty}\binom{k+j}{j}t^j(1-t)^k n_{+}(1)f_{\widehat{Q}}(k+j).
\end{equation}
The fitted density $f_{\widehat{Q}}(x)$
is used to estimate $f_Q(x)$ and yield $\hat{\phi}_j(t)$ while
the empirical density $\hat{f}_Q(x)=n_x(1)/n_+(1)$ is used to estimate
$f_Q(x)$ and yield $\tilde{\phi}_j(t)$.

The function $G(\phi(t))$ can be estimated by
$G(\tilde{\phi}(t))$ and  $G(\hat{\phi}(t))$.
The estimator $G(n(1))$ is reproduced
by $G(\tilde{\phi}(1))=G(n(1))$.
A bootstrap procedure is recommended for construction of
confidence intervals for $G(\phi(t))$: 
sampling $n_+^{\star}(1)$ from its estimated binomial density
and sampling $n^{\star}(1)$
from $L(\widehat{Q}, n_+^{\star}(1))$.
A lower confidence limit for $G(\phi(t))$
is also a lower confidence limit for  $s$ when $G(\phi(t))$
is a lower bound to $s$, e.g.,  $\phi_+(t) $ and  $G_c(\phi(t))$.

It is difficult to estimate $\phi_1(t)$ reliably  when $t$ is relatively large.
One reason is that, although $\gamma_u>0$ in $Q$ for all $u$, 
the smallest support point (say $\hat{\gamma}_1$)
of  $\widehat{Q}$ might be close or identical to zero.
When $\hat{\gamma}_1=0$, it is easily shown that
\[
\theta_j(t,\widehat{Q})=I(j=1)\hat{\omega}_1t+
\sum_{u=2}^{\hat{\nu}} 
\hat{\omega}_u\frac{\exp(-\hat{\gamma}_u t)(\hat{\gamma}_ut)^j}
{\{1-\exp(-\hat{\gamma}_u)\}j!}.
\] 
When $t$ is sufficiently
large,  $\hat{\phi}_1(t)$ will increase approximately linearly 
but each $\hat{\phi}_j(t)$ with $j\geq 2$ will approach zero.
This fact explains the observation 
that $\hat{\phi}_+(t)$ is approximately linear 
for a large $t$  (Mao 2005). The estimator 
$G(\hat{\phi}(t))$ 
might also be driven up to infinity as $t$ increases.
For example, if $\hat{\gamma}_1=0$, then 
there is $\beta\geq 2$ with $\hat{\gamma}_{\beta}<\hat{\gamma}_u$ 
for all $u\geq 2$ and $u\neq \beta$, and
\[
\lim_{t\to\infty}\frac{G_c(\hat{\phi}(t))}
{\exp(\hat{\gamma}_{\beta} t)}=
\frac{n_+(1)\hat{\omega}_1^2\{1-\exp(-\hat{\gamma}_{\beta})\}}
{2\hat{\omega}_{\beta} \hat{\gamma}_{\beta}^2},
\]
i.e.,  $G_c(\hat{\phi}(t))$ increases approximately
exponentially for a large $t$. However, our likelihood-based method
can be useful for relatively small $t$ (e.g., $t\in[1,3]$ with $t_0=1$,
the range of $t$ that serves practical purposes).

Finally we turn to the multinomial model. 
Let $X_i(h)$ be the number of individuals from  species $i$ in a sample of size $h$
and $m_j(h)=\sum_{i=1}^sI(X_i(h)=j)$.
This means that $X_i(h(t))=Y_i(t)$ and $m_j(h(t))=n_j(t)$.
Note that
\[
E\{m_j(h)\}=\sum_{i=1}^s\binom{h}{j}p_i^j(1-p_i)^{h-j}.
\]

Let $a=h(1)$ be the number of sampled individuals during $[0,1]$.
For $h=1$, 2, \dots, $a$, one has
\begin{equation}\label{eqn:rhohat}
\hat{m}_j(h)=\sum_{k=0}^{a-h}\binom{h}{j}
\binom{a-h}{k}\binom{a}{k+j}^{-1}m_{k+j}(a), j=1, 2,\dots, h,
\end{equation}
which is based on the following identity (Good and Toulmin 1956)
\begin{equation}\label{eqn:id2}
E\{m_j(h)\}
=\sum_{k=0}^{a-h}\binom{h}{j}\binom{a-h}{k}\binom{a}{k+j}^{-1}E\{m_{k+j}(a)\},
j=1, 2,\dots, h.
\end{equation}

In the ecology literature, a randomization procedure is usually used.
It is an approximation to an enumeration procedure: taking all subsamples of size $h$,
calculate $m_j(h)$ with $j\geq h$ for each subsample 
and obtain their $\bar{m}_j(h)$. 

\begin{proposition}For $j=1$, 2, \dots, $h$ and $h=1$, 2, \dots, $a$,
\begin{equation}\label{eqn:rhobar}
\bar{m}_j(h)=\sum_{k=0}^{a-h}\binom{k+j}{j}\binom{a-k-j}{h-j}\binom{a}{h}^{-1}m_{k+j}(a).
\end{equation}
\end{proposition}

Hurlbert (1971) found the analytic expression of
$\bar{m}_+(h)=\sum_{j=1}^h\bar{m}_j(h)$,
\[
\bar{m}_+(h)=\sum_{x=1}^{a}m_{x}(1)
-\sum_{x=1}^{a-h}\binom{a-h}{x}\binom{a}{x}^{-1}m_{x}(1).
\]

By comparing (\ref{eqn:rhohat}) and (\ref{eqn:rhobar}), 
it is clear that $\bar{m}_j(h)=\hat{m}_j(h)$  because
\[
\frac{\binom{k+j}{j}\binom{a-k-j}{h-j}}{\binom{a}{h}}
=\frac{\binom{h}{j}\binom{a-h}{k}}{\binom{a}{k+j}}
=\frac{(k+j)!(a-k-j)!h!(a-h)!}{j!k!(h-j)!(a-k-h)!a!}.
\]

Although the identity in (\ref{eqn:GTsource}) holds
for all $t>0$, the identity in (\ref{eqn:id2})
does not hold for $h>a$. One can obtain
an approximation to $E\{m_j(h)\}$ as a function
of those $E\{m_j(a)\}$ and develop a biased estimator for $E\{m_j(h)\}$.

Since $m_{b}(a)=n_b(1)$, we can write $\bar{m}_j(h)$ as
\begin{equation}\label{eqn:rhobar2}
\bar{m}_j(h)=\sum_{b=j}^{\min(a-h+j,d)}
\binom{b}{j}\binom{a-b}{h-j}\binom{a}{h}^{-1}n_{b}(1).
\end{equation}
The number of sampled individuals during $[0,h/a]$ is about $h$.
We consider comparing the estimators $\tilde{\phi}_j(h/a)$ in (\ref{eqn:rep1})
and $\bar{m}_j(h)$  in (\ref{eqn:rhobar2}).
Clearly $\bar{m}_j(h)=\tilde{\phi}_j(h/a)=0$ when $j>d$. When $j\leq d$,
write $\tilde{\phi}_j(h/a)-\bar{m}_j(h)=\epsilon_1+\epsilon_2$, where
\begin{gather*}
\epsilon_1=\sum_{b=\min(a-h+j,d)+1}^d\binom{b}{j}(h/a)^j\{(a-h)/a\}^{b-j}n_b(1),\\
\epsilon_2=\sum_{b=j}^{\min(a-h+j,d)}\binom{b}{j}
\left[(h/a)^j\{(a-h)/a\}^{b-j}-
\prod_{u=0}^{j-1}\frac{h-u}{a-u}\prod_{w=0}^{b-j-1}\frac{a-h-w}{a-j-w}
\right]n_b(1).
\end{gather*}
Note that $\epsilon_1=0$ when $a-h+j\geq d$. When  $a-h+j<d$,
$h/a$ is close to one because $d\lll a$, which implies that $\epsilon_1\approx 0$.
By simple algebra, one can also find that $\epsilon_2\approx 0$.
Conclude that $\bar{m}_j(h)\approx\tilde{\phi}_j(h/a)$.
When the $\bar{m}_j(h(t))$ are used to estimate
$G(\phi(t))$, the resulting estimator will be
close to $G(\tilde{\phi}(t))$.
For example, $\bar{m}_+(h)$ and $\tilde{\phi}_+(h/a)$
are close to one another (Brewer and Williamson 1994).

\section{A real example}
We consider a real example
from Miller and Wiegert (1989) that 
concerns plant species in the central
Appalachian region. This example was also
investigated in Shen et al. (2003).
There were $n_+(1)=188$ species
identified from $a=h(1)=1008$ individuals
with $n_x(1)=61$, 35,  18,  12,  15,  4,   8,   4,  
5,   5,   1,   2,   1,   2,   3,  
2,   1,   2,   1,   1,   1,   1,   1,   1 and  1
at $x=1$, 
  2,   3,   4,   5,   6,   7,   8,   9,   10,  11,  12,  13,  14,  15, 
 16,  19,  20,  22,  29,  32,  40,  43,  48 and  67.

\begin{figure}[htbp]
  \centering
\includegraphics[angle=270,width=10cm]{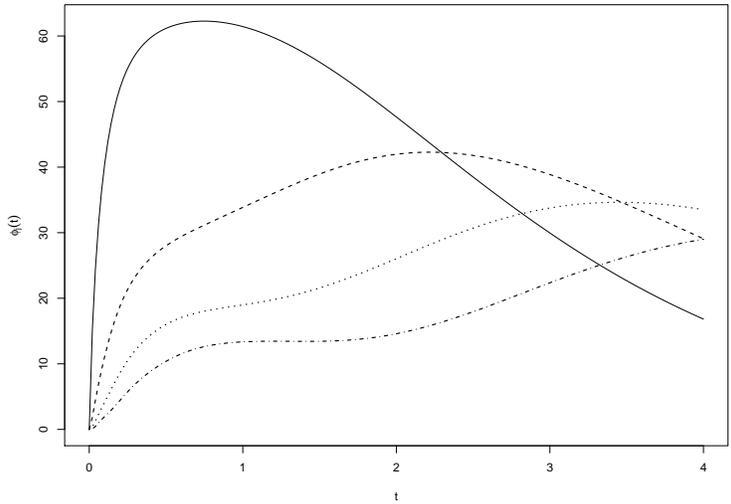}
\caption{The likelihood-based estimates $\hat{\phi}_j(t)$
of the expected counts $\phi_j(t)$ for 
$j=1$ (solid), 2 (dashed), 3 (dotted) and 4 (dot-dashed).}
  \label{fig:ec}
\end{figure}

\begin{figure}[htbp]
  \centering
\includegraphics[angle=270,width=10cm]{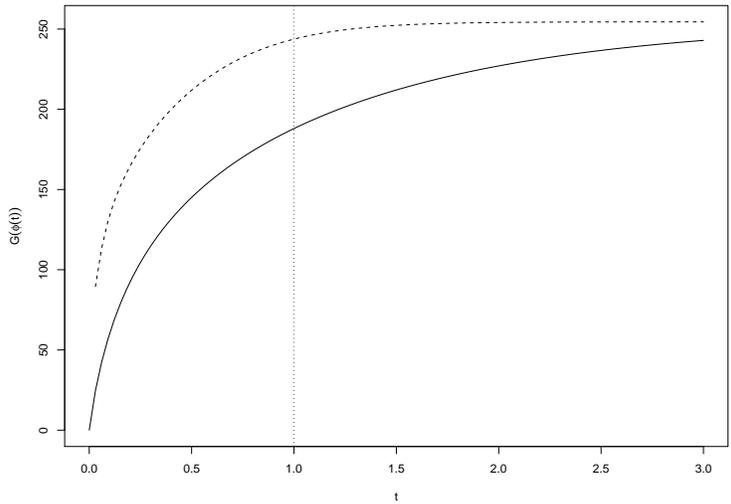}
\caption{The likelihood-based estimates $G_c(\hat{\phi}(t))$ 
(dashed) and $\phi_+(t)$ (solid).}
  \label{fig:gsac}
\end{figure}

The NPMLE $\widehat{Q}$ is shown in
Table~\ref{tab:plant-Q}. 
The estimates
$\hat{\phi}_j(t)$, $\phi_+(t)$ and
$G_c(\hat{\phi}(t))$ are shown 
in Figures~\ref{fig:ec} and \ref{fig:gsac}.
We also compare  $\tilde{\phi}_j(h/a)$
and $\bar{m}_j(h)$ for $1\leq h\leq a$,
and $\tilde{\phi}_j(t)$
and $\hat{\phi}_j(t)$ for $0\leq t\leq 1$.
The results are shown in Table~\ref{tab:comparison}.
We also calculate
$\max_{0\leq t\leq 1}|\tilde{\phi}_+(t)-\hat{\phi}_+(t) |=0.06$
and $\max_{0\leq t\leq 1}|G_c(\tilde{\phi}(t))
-G_c(\hat{\phi}(t))|=2.58$.
Note that $\tilde{\phi}_j(h/a)$
and $\bar{m}_j(h)$ have little difference.
The difference between $\tilde{\phi}_j(t)$
and $\hat{\phi}_j(t)$  comes from the difference
between  $n_x(1)$ and $n_+(1)f_{\widehat{Q}}(x)$,
e.g.,  $n_x(5)=15$ and $n_+(1)f_{\widehat{Q}}(5)=10.4$.
Although $\tilde{\phi}_j(t)$
can be computed for $t>1$,
it becomes inadmissible even for some $t<2$,
e.g., $\tilde{\phi}_2(1.57)=-2.6$
and $\tilde{\phi}_4(1.57)=-192.5$, 
$G_c(\tilde{\phi}(1.57))=-497.1$.
To construct lower confidence limits for
$G_c(\phi(t))$, we generate 400 bootstrap
resamples. For example, the bootstrap 95\% lower confidence
limits for $G_c(\phi(t))$ at $t=1$, 1.2, 1.4, 1.6, 1.8 and 2.0
are 218.1, 220.6, 221.5, 222.2, 222.4 and 222.4
while the estimates $G_c(\hat{\phi}(t))$
are  243.7, 248.7, 251.5, 252.9, 253.7 and 254.1. 
Note that an upper confidence limit at a relatively
large $t$ is usually noninformative.
For example, the 95\%  upper confidence
limits for $G_c(\phi(t))$ at $t=2$ and $t=3$
are 811.8 and 2230.8 respectively,
much larger than the corresponding
lower confidence limits 222.4  and 222.6.

In order to evaluate the likelihood-based method,
we consider simulation under various
combinations of $Q$ and $s$. 
We find that
the distribution of $\hat{\phi}_j(t)$
is right skewed when $t>2$
and in particular,
the distribution of $\hat{\phi}_1(t)$
has a long right tail for a large $t$,
like $\hat{\phi}_+(t)$ and  $G_c(\hat{\phi}(t))$
although the 3rd quartile
of $G_c(\hat{\phi}(t))$ increases
faster than
that of $\hat{\phi}_+(t)$ or $\hat{\phi}_1(t)$.
In the future, we will consider 
generalized SACs for various
nonparametric estimators 
(e.g., Chao and Bunge 2002).

\begin{table}[ht]
\begin{center}
\caption{The NPMLE $\widehat{Q}$ with $\hat{\nu}=7$ from the plant data. }
\label{tab:plant-Q}
\begin{tabular}{rrrrrrrr}
\hline
$\hat{\gamma}_u$ & 0.864 & 3.554 & 7.412 & 15.306 & 30.564 & 41.892 & 66.416 \\
$\hat{\omega}_u$ & 0.475 & 0.260 & 0.158 & 0.074 & 0.010 & 0.017 & 0.005 \\
\hline
\end{tabular}
\end{center}
\end{table}

\begin{table}[ht]
\begin{center}
\caption{Comparison of three types of estimates $\tilde{\phi}_j(t)$,
$\hat{\phi}_j(t)$ and $\bar{m}_j(h)$
with $\Delta_{j}=\max_{1\leq h\leq a}|\tilde{\phi}_j(h/a)-\bar{m}_j(h)|$
and $D_{j}=\max_{0\leq t\leq 1}|\tilde{\phi}_j(t)-\hat{\phi}_j(t)|$.}
\label{tab:comparison}
\begin{tabular}{rrrrrrrrrrr}
\hline
$j$ & 1 & 2 & 3 & 4 & 5 & 6 & 7 & 8 & 9 & 10 \\
\hline
$\Delta_{j}$
 & 0.13 & 0.03 & 0.02 & 0.01 & 0.01 & 0.01 & 0.01 & 0.01 & 0.01 & 0.01 \\
$D_{j}$
 & 0.44 & 1.14 & 1.01 & 1.35 & 4.56 & 4.26 & 1.44 & 1.21 & 1.03 & 1.78 \\
\hline
\end{tabular}
\end{center}
\end{table}

\nocite{*}
\bibliographystyle{apalike}
\bibliography{ref}

\section*{Appendix}

To prove Proposition 1, write
\begin{align*}
\frac{\phi_j(t)}{\phi_+(1)}
&=\frac{s\sum_{u=1}^{\nu}\pi_u \exp(-\gamma_u t)(\gamma_u t)^j(j!)^{-1}}
{s-s\sum_{w=1}^{\nu}\pi_w \exp(-\gamma_w t)}\\
&=\sum_{u=1}^{\nu}\frac{\pi_u \{1-\exp(-\gamma_u t)\}}
{\sum_{w=1}^{\nu}\pi_w\{1- \exp(-\gamma_w t)\}}\cdot
\frac{\exp(-\gamma_u t)(\gamma_u t)^j}{\{1-\exp(-\gamma_u t)\}j!}.
\end{align*}

To prove Proposition 2, let the individuals be labeled by $j=1$, 2, \dots, $a$
and $Z_{ij}=I(\text{individual } j \text{ is from species } i)$.
A subsample $\omega$ consists of $h$ individuals. 
Let $\Omega$ be the set of all such subsamples.
With $\binom{\alpha}{\beta}=0$ if $\alpha<\beta$, write
\begin{multline*}
\tbinom{a}{h}\bar{n}_j(h)=\sum_{\omega\in\Omega}
\sum_{i=1}^sI\Bigl(\sum_{r\in\omega}Z_{ir}=j\Bigr)
=\sum_{t=0}^a\sum_{\{i:Y_i(a)=t\}}
\sum_{\omega\in\Omega}I\Bigl(\sum_{r\in\omega}Z_{ir}=j\Bigr)\\
=\sum_{t=0}^a\sum_{\{i:Y_i(a)=t\}}\tbinom{t}{j}\tbinom{a-t}{h-j}
=\sum_{t=j}^{a-h+j}\tbinom{t}{j}\tbinom{a-t}{h-j}n_t(a)
=\sum_{k=0}^{a-h}\tbinom{k+j}{j}\tbinom{a-k-j}{h-j}n_{k+j}(a).
\end{multline*}

\end{document}